\documentclass[a4paper,12pt,twoside]{article}
\usepackage{amsbsy}
\usepackage{amssymb}
\usepackage{amsfonts}
\usepackage{amsmath}
\usepackage{amsopn}
\usepackage{amsthm}
\usepackage[english]{babel}
\usepackage[T1]{fontenc}

\newtheorem{teo}{Theorem}
\newtheorem{fact}{Fact}
\newtheorem{lemma}{Lemma}
\newtheorem{prop}{Proposition}

\theoremstyle{definition}
\newtheorem{remark}{Remark}
\theoremstyle{definition}

 \DeclareMathOperator{\er}{\mathbf{E}}
\DeclareMathOperator{\pr}{\mathbf{P}}
\DeclareMathOperator{\erc}{\mathbf{E}^{\scriptscriptstyle{\downarrow}}}
\DeclareMathOperator{\prc}{\mathbf{P}^{\scriptscriptstyle{\downarrow}}}

\DeclareMathOperator{\prcp}{\sideset{}{^{\searrow}}\pr}

\DeclareMathOperator{\re}{\mathbb{R}}

\DeclareMathOperator{\ke}{\mathrm{q}}
\DeclareMathOperator{\p}{\mbox{\rm
I\hspace{-0.02in}P}} \DeclareMathOperator{\e}{\mbox{\rm
I\hspace{-0.02in}E}}

\DeclareMathOperator{\pc}{\mbox{\rm
I\hspace{-0.02in}P}^{\scriptscriptstyle{\downarrow}}}
\DeclareMathOperator{\pcx}{\mbox{\rm
I\hspace{-0.02in}P}^{\scriptscriptstyle{\downarrow}}_{\mathnormal
x}}

\newcommand{\R}{\mbox{\rm I\hspace{-0.02in}R}}
\newcommand{\ep}{\mbox{\raisebox{-1.1ex}{$\rightarrow$}\hspace{-.12in}$\e$}}
\newcommand{\poxi}{\mbox{\raisebox{-1.1ex}{$\rightarrow$}\hspace{-.12in}$\xi$}}
\newcommand{\taup}{\mbox{\raisebox{-1.1ex}{$\rightarrow$}\hspace{-.12in}$\tau$}}
\newcommand{\ppx}{\mbox{\raisebox{-1.1ex}{$\rightarrow$}\hspace{-.12in}$\p^x$}}

\newcommand{\ppxy}{\mbox{\raisebox{-1.1ex}{$\rightarrow$}\hspace{-.12in}$\p^x_y$}}
\newcommand{\pp}{\mbox{\raisebox{-1.1ex}{$\rightarrow$}\hspace{-.12in}$\p$}}
\newcommand{\ed}{\stackrel{(d)}{=}}
\newcommand{\eqdef}{\stackrel{\mbox{\tiny(def)}}{=}}

\numberwithin{equation}{section}

\title{\textbf{On some transformations between positive self--similar Markov processes}}
\author{\textbf{Lo\"\i c CHAUMONT}\thanks{Laboratoire de Probabilit\'es et
Mod\`eles Al\'eatoires, Universit\'e Pierre et Marie Curie, 4, Place
Jussieu - 75252 {\sc Paris Cedex 05, France.}} \and  \textbf{V\'\i
ctor RIVERO}\thanks{CIMAT A.C.
 Calle Valenciana s/n, C.P.36240, Guanajuato, Gto. Mexico; and Universit\'e Paris X,
 Nanterre, MODAL'X, 200, Avenue
de la R\'epublique, 92001, Nanterre Cedex, France; Email:
rivero@cimat.mx}}
\date{\footnotesize This version: \today}

\begin{document}

\maketitle
\begin{abstract}
A path decomposition at the infimum for positive self-similar Markov
processes (pssMp) is obtained. Next, several aspects of the
conditioning to hit 0 of a pssMp are studied. Associated to a given
a pssMp $X,$ that never hits 0, we construct a pssMp
$X^{\downarrow}$ that hits 0 in a finite time. The latter can be
viewed as $X$ conditioned to hit 0 in a finite time and we prove
that this conditioning is determined by the pre-minimum part of $X.$
Finally, we provide a method for conditioning a pssMp that hits 0 by
a jump to do it continuously.
\end{abstract}


\noindent
\noindent\textbf{MSC:} 60 G 18 (60 G 17).
\begin{section}{Introduction}

This work concerns positive self--similar Markov processes (pssMp),
that is $[0,\infty[$-valued strong Markov processes that have the
scaling property: there exists an $\alpha>0$ such that for any
$0<c<\infty$, $$\left\{ \left(cX_{tc^{-1/\alpha}}, t\geq
0\right),\p_x\right\}\ed\left\{ \left( X_{t}, t\geq
0\right),\p_{cx}\right\},\qquad x>0.$$ This class of processes has
been introduced by Lamperti~\cite{la} and since then studied by
several authors, see
e.g.~\cite{bc,by,by2,ch3,chaumontcaballero,ri2,ri}. We will make
systematic use of a result due to Lamperti that establishes that any
pssMp is the exponential of a L\'evy process time changed, this will
be recalled at Section~\ref{lamptrasform}.

Some of the motivations of this work are some path decompositions
and conditionings that can be deduced from ~\cite{ch2,ch3} and that
we will recall below, in the particular case where the positive
self--similar Markov process is a stable L\'evy process conditioned
to stay positive.

Let $\widetilde{\p}$ be a law on the the  space of c\`adl\`ag paths
under which the canonical process $X,$ is an $\alpha$-stable L\'evy
process, $0<\alpha\leq 2$, i.e. a process with independent and
stationary increments that is $1/\alpha$-self--similar. Associated
to this process we can construct a pssMp, say $(X,\p),$ that can be
viewed as $(X,\widetilde{\p})$ conditioned to stay positive. The
construction can be performed either via the
Tanaka--Doney~\cite{doneytanaka} path transform of
$(X,\widetilde{\p})$ or as an $h$-transform of the law of
$(X,\widetilde{\p})$ as in \cite{ch2,chaumontdoney} or also, in the
spectrally one sided case, via Bertoin's
transformation~\cite{be:splitting}.

Another interesting process related to $(X,\widetilde{\p})$ is
$(X,\sideset{}{_{\cdot}^{\downarrow}}\p),$  which was introduced
in~\cite{ch2}, can be viewed as $(X,\widetilde{\p})$ conditioned to
hit 0 continuously and is constructed via an h-transform of
$(X,\widetilde{\p})$ killed at its first hitting time of
$]-\infty,0].$

Using the results of Millar~\cite{mi2}, in \cite{ch2} it has been
proved the following results for $(X,\p),$ relating $(X,\pc)$ and
$(X,\p)$ started at 0, with the pre and post minimum parts of
$(X,\p).$

\begin{fact}\label{stablepathdec}
Let $I^{X}=\inf\{X_s, s>0\}$ and $m=\sup\{t>0 : X_{t-}\wedge
X_t=I^{X}\}.$ Under $\p$, the pre-minimum part of $X$,
i.e.~$\{X_t,0\leq t<m\},$ and the post minimum part of $X$,
i.e.~$\{X_{m+t}, t>0\}$ are conditionally independent given the
value of $I^X.$ For any $x>0,$ under $\p_x,$ conditionally on
$I^{X}=y,$ $0<y\leq x,$ the law of the former is
$(X+y,\sideset{}{_{x-y}^{\downarrow}}\p)$ and that of the later is
$(X+y,\sideset{}{_{0+}}\p),$ where $\sideset{}{_{0+}}\p$ is the
limiting law of $(X,\p_{\cdot})$ as the starting point tends to
$0$, $\p_x\xrightarrow{w} \sideset{}{_{0+}}\p$ as $x\to 0+.$
\end{fact}
Furthermore, it can be verified using the previous result, and it
is intuitively clear, that under $\p_x$ the law of the pre-minimum
(respectively, post-minimum) of $X,$ conditionally on the event
$\{I^{X}<\epsilon\},$ converges as $\epsilon\to 0,$ to the law
$\sideset{}{_x^\downarrow}\p,$ respectively $\sideset{}{_{0+}}\p,$
in the sense that,
\begin{fact}\label{asymdecomp}$\displaystyle
\lim_{\epsilon\to 0+}\p_x(F\cap\{t<m\}, G\circ\theta_m | I^X <
\epsilon)=\sideset{}{_x^{\downarrow}}\p(F\cap\{t<T_0\})\sideset{}{_{0+}}\p(G),\quad
F\in \mathcal{G}_{t}, G\in\mathcal{G}_{\infty},$ where
$\{\mathcal{G}_t,t\geq 0\}$ is the natural filtration generated by
$X.$
\end{fact}
Our first purpose is to extend
Facts~\ref{stablepathdec}~\&~\ref{asymdecomp}, to a larger class of
positive self-similar Markov processes. That is the content of
sections~\ref{pathdecomp} \& \ref{prepostas}, respectively.

Here is another interpretation of the law
$\sideset{}{^\downarrow}\p$. Let $\p^{0}$ be the law of  the
process $(X,\widetilde{\p})$ killed at its first hitting time of
$]-\infty,0].$ The process $(X,\p^{0})$ still has the strong
Markov property and inherits the scaling property from
$(X,\widetilde{\p}),$ so it is a pssMp and it hits 0 in a finite
time. Moreover, whenever $(X,\widetilde{\p})$ has negative jumps,
the process $(Y,\p^{0})$ hits $0$ for the first time with a
negative jump:
$$\p^{0}_x(T_0<\infty, X_{T_0-}>0)=1,\qquad \forall x>0,$$ where
$T_0=\inf\{t>0: X_{t}=0\}.$ It has been proved in \cite{ch2}, that
$\pc$ is an $h$ transform of $\p^{0}$ via the excessive function
$x\mapsto x^{\alpha(1-\rho)-1},$ $x>0,$ where $\rho$ is the
positivity parameter of  $(X,\widetilde{\p}),$
$\rho=\widetilde{\p}(X_1\geq 0).$ Furthermore,  $(X,\pc)$ hits 0
continuously and in a finite time, i.e.:
$$\sideset{}{^{\downarrow}_x}\p(T_0<\infty, X_{T_0-}=0)=1,\qquad \forall x>0,$$
and Proposition~3 in~\cite{ch2} describes a relationship between
$\pc$ and $\p^{0}$ that allows us to refer to $\pc$ as the law of
$(X,\p^{0})$ conditioned to hit $0$ continuously. The latter
conditioning is performed by approximating the set
$\{I^{Y^{0}}=0\}$ by the sequence of sets $\{I^{Y^{0}}<\epsilon\},
\epsilon >0.$

In Section~\ref{condpssMphito0cont}, we obtain an analogous result
for a larger class of self--similar Markov processes. Namely
those associated to a L\'evy process killed at an independent
exponential time and which satisfy a Cram\'er's type condition.
Furthermore, an alternative method for conditioning a self-similar
Markov process that hits 0 by a jump, to hit 0 continuously, is
provided by making tend to 0 the height of the jump by which the
process hits the state 0.

The approach used to aboard these problems is based on Lamperti's
representation between real valued L\'evy processes and pssMp
which we recall in the following section.
\end{section}

\begin{section}{Some preliminaries on pssMp}\label{lamptrasform}
Let $\mathbb{D}$ be the space of c\`adl\`ag paths defined on
$[0,\infty)$, with values in $\R\cup\Delta$, where $\Delta$ is a
cemetery point. Each path $\omega\in\mathbb{D}$ is such that
$\omega_t=\Delta$, for any
$t\ge\inf\{t:\omega_t=\Delta\}:=\zeta(\omega)$. As usual we extend
the functions $f:\re\to\re$ to $\re\cup\Delta$ by $f(\Delta)=0.$ The
space $\mathbb{D}$ is endowed with the Skohorod topology and its
Borel $\sigma$-field. We will denote by $X$ the canonical process of
the coordinates and $({\cal F}_t)$ will be the natural filtration
generated by $X$. Moreover, let $\pr$ be a reference probability
measure on $\mathbb{D}$ under which the process, $\xi,$ is a L\'evy
process; we will denote by $(\mathcal{D}_t,t\geq 0),$ the complete
filtration generated by $\xi.$

Fix $\alpha> 0$ and let $(\p_x,x>0)$ be the laws of an
$\alpha$-pssMp associated to $(\xi,\pr)$ via the Lamperti
representation. Formally, define
$$A_t=\int^t_0\exp\{(1/\alpha)\xi_s\}ds,\qquad t\geq 0,$$ and let
$\tau(t)$ be its inverse, $$\tau(t)=\inf\{s>0 : A_s>t\},$$ with the
usual convention, $\inf\{\emptyset\}=\infty.$ For $x>0,$ we denote
by $\p_x$ the law of the process
$$x\exp\{\xi_{\tau(tx^{-1/\alpha})}\}, \qquad t>0,$$ with the
convention that the above quantity is $\Delta$ if
$\tau(tx^{-1/\alpha})=\infty.$ The Lamperti representation ensures
that the laws $(\p_x, x>0)$ are those of a pssMp with index of
self-similarity $\alpha$.

Besides, recall that any L\'evy process $(\xi,\pr)$ with lifetime
has the same law as a L\'evy process with infinite lifetime that has
been killed at a rate $\ke\geq 0.$ It follows that $T_0=\inf\{t>0:
X_{t}=0\}$ has the same law under $\p_x$ as $x^{1/\alpha}A_{\zeta}$
under $\pr$ with
\begin{equation*}\label{eq:1ch4}
A_{\zeta}=\int^{\zeta}_0\exp\{(1/\alpha)\xi_s\}ds.
\end{equation*}
So, if $\ke>0,$ then the random variable $A_{\zeta}$ is a.s. finite;
while in the case $\ke=0,$ we have two possibilities, either
$A_{\zeta}$ is finite a.s. or infinite a.s.; the former happens if
and only if $\lim_{t\to\infty}\xi_t=-\infty,$ a.s. and the latter if
and only if $\limsup_{t\to\infty}\xi_t=\infty,$ a.s.

Lamperti proved that any pssMp can be constructed this way and
obtained the following classification of pssMp's:
\begin{itemize}
\item[(LC1)] $\ke>0,$ if and only if
\begin{equation}\label{LC1}\p_x(T_0<\infty,\ X_{T_0-}>0,\ X_{T_0+s}=0,\ \forall
s \geq 0)=1,\qquad \text{for all}\ x>0.\end{equation}
\item[(LC2)] $\ke=0$ and $\lim_{t\to\infty}\xi_t=-\infty$ a.s. if and
only if
\begin{equation}\label{LC2}\p_x(T_0<\infty,\ X_{T_0-}=0,\ X_{T_0+s}=0,\ \forall
s \geq 0)=1,\qquad \text{for all}\ x>0,\end{equation}
\item[(LC3)] $\ke=0$ and $\limsup_{t\to\infty}\xi_t=\infty$ a.s.
if and only if
\begin{equation}\label{LC3}\p_x(T_0=\infty)=1,\qquad \text{for all}\ x>0.
\end{equation}
\end{itemize}

Observe that without loss of generality we can and we will suppose
that $\alpha=1$ in Lamperti's construction of pssMp because all our
results may trivially be extended to any $\alpha>0$ by considering
$X^\alpha$ which is a pssMp with index of self-similarity $\alpha$.

In this work we will be mostly interested by those pssMp that
belong to the class LC3; nevertheless, in
Section~\ref{condpssMphito0cont} we will prove that some elements
of the class~LC1 can be transformed into elements of the class
LC2.
\end{section}

\begin{section}{Path decomposition at the minimum}\label{pathdecomp}
\setcounter{equation}{0} \noindent We suppose throughout this
section that $(\xi,\pr)$ is a L\'evy process with infinite
lifetime which drifts to $+\infty$, that is
$\lim_{t\rightarrow+\infty}\xi_t=+\infty$, a.s. We start by
recalling a William's type path decomposition of $(\xi,\pr)$ at
its minimum. Let $I^\xi=\inf_{t\geq0}\xi_t$ and
$\rho=\sup\{t:\xi_t\wedge \xi_{t-}=I^\xi\}$. We define the post
minimum process  as
\[\poxi\eqdef(\xi_{\rho+t}-I^\xi,\,t\ge0)\,.\]
The following result is due to Millar \cite{mi2}, proposition 3.1
and Theorem 3.2. \noindent
\begin{teo}\label{decxi}
The pre-minimum process $((\xi_t,\,t< m),\pr)$ and the
post-minimum process $(\poxi\,,\pr)$ are independent. Moreover,
the three following exhaustive cases hold:
\begin{itemize}
\item[$(i)$] $0$ is regular for both $(-\infty,0)$ and
$(0,\infty)$ and $\pr$-a.s., there is no jump at the minimum,

\item[$(ii)$] $0$ is regular for $(-\infty,0)$ but not for
$(0,\infty)$ and $I^\xi=\xi_{\rho-}<\xi_\rho$ , $\pr$-a.s.

\item[$(iii)$]  $0$ is regular for $(0,\infty)$ but not for
$(-\infty,0)$ and $I^\xi=\xi_\rho<\xi_{\rho-}$, $\pr$-a.s.
\end{itemize}
In any case under $\pr$, the process $(\xi_t,\,t< m)$ and $\poxi$
are also conditionally independent given $I^\xi$ and the process
$\poxi$ is strongly Markovian.
\end{teo}
\noindent Actually, Millar's result is much more general and
asserts that for any Markov process, which admits a minimum, the
pre-minimum process and the post minimum process are conditionally
independent given both the value at the minimum and the subsequent
jump and the post-minimum process is strongly Markovian. When this
Markov process is a pssMp that belongs to the class (LC3), we may
complete Millar's result as in the following proposition. First of
all, observe that $X$ derives towards $+\infty$ as well as $\xi,$
and so the following are well defined $I^X=\inf_{t\geq 0}X_t$ and
$m=\sup\{t:X_t\wedge X_{t-}=I^X\}$.
\begin{prop}\label{pathdecpssmp} For any $x>0,$ under $\p_x,$ the processes
$(X_t,t<m)$ and $(X_{t+m},\,t\ge0)$
are conditionally independent given $I^X$, and with the
representation given by Lamperti's
transformation~(Section~\ref{lamptrasform}), we have
\begin{eqnarray}\left((X_t,0\le t<m),\p_x\right)&=&
\left(\left(x\exp\xi_{\tau(t/x)},\,0\le t<x\int^\rho_0
\exp{\xi_s}\mathrm{d}s\right),\pr\right),\label{dec1}\\
\left((X_{t+m},\,t\ge0),\p_x\right)&=&\left(\left(I^X\exp
\poxi_{\;\taup\,(t/I^X)},\,t\ge0\right),\pr\right)\,,\label{dec2}
\end{eqnarray}
where
$\taup(t)=\inf\left\{s:\int_0^s\exp\left(\poxi_{\;u}\right)\,du>t\right\}$,
for $t\ge0$.
\end{prop}

\begin{proof} The expression of the pre-minimum part of
$(X,\p_x)$ follows directly from Lamperti's transformation
(Section~\ref{lamptrasform}). Note that in particular, since $\tau$
is a continuous and strictly increasing function, one has:
\begin{equation}\label{time}
A_{\rho}=\int^\rho_0\exp{\xi_s}\mathrm{d}s\,,\;\;\;\tau(A_\rho)=
\rho\,,\;\;\;xA_\rho=m\;\;\;\mbox{and}\;\;\;I^X=x\exp
I^\xi\,.
\end{equation}
To express the post-minimum part of $(X,\p_x)$, first note that
\[X_{m+t}=x\exp\xi_{\tau(A_\rho+t/x)}\,,\;\;t\ge0\,.\]
Then we can write the time change as follows:
\begin{eqnarray*}
\tau(A_\rho+t/x)&=&\inf\{s>0:\int_0^s\exp\xi_u\,du>A_\rho+t/x\}\\
&=&\inf\{s>\rho:\int_0^{s-\rho}\exp\xi_{u+\rho}\,du>t/x\}\\
&=&\inf\{s>0:\int_0^s\exp\poxi_{\,u}\,du>(t/x)\exp(-I^\xi)\}+\rho\\
&=&\taup((t/x)\exp(-I^\xi))+\rho=\taup(t/I^X)+\rho\,,
\end{eqnarray*}
so that
\[\xi_{\tau(A_\rho+t/x)}=\xi_{\taup(t/I^X)+\rho}=\poxi_{\taup(t/I^X)}+I^\xi\]
and the expression (\ref{dec2}) for the post-minimum part of
$(X,\p_x)$ follows.

From (\ref{dec1}), we see that $(X_t\,,t<m)$ is a measurable
functional of $(\xi_t\,,t<\rho)$ and from (\ref{dec2}),
$(X_{m+t},\,t\ge0)$ is a functional of $I^X$ and $\poxi\;$. Since
$I^X=x\exp I^\xi$, the conditional independence follows from
Theorem \ref{decxi}.\end{proof}

\noindent When $X$ has no positive jumps (or equivalently when
$\xi$ has no positive jumps), it makes sense to define the last
passage time at level $y\ge x$ as follows
\[\sigma_y=\sup\{t:X_t=y\}\,.\]
Then the post-minimum process of $X$ becomes more explicit as the
following result shows; its proof is an easy consequence of
Proposition~\ref{pathdecpssmp}.

\begin{prop}
Let $y\leq x$. Conditionally on $I^X=y$, the post-minimum process
$(X_{t+m},\,t\ge0)$ has the same law as $(X_{\sigma_y+t},\,t\ge0)$,
and
\begin{equation}\label{postsig}
((X_{\sigma_y+t},\,t\ge0),\p_x)\ed
\left(\left(y\exp\poxi_{\;\taup(t/y)},\,t\ge0\right),\pr\right)\,.
\end{equation}
\end{prop}

As we have just seen, the post-minimum process of $(X,\p)$ can be
completely described using the underlying L\'evy process
$(\xi,\pr)$ conditioned to stay positive $(\poxi,\pr)$.
Nevertheless, the description of the pre-minimum obtained
in~(\ref{dec1}) is not so explicit. So our next purpose is to make
some contributions to the understanding of the pre-minimum process
of a positive self-similar Markov process.

Let us start by the case where the process $(X,\p)$ (or equivalently
the underlying L\'evy process $(\xi,\pr)$) has no negative jumps,
because in this case we can provide a more precise description of
the pre-minimum process using known results for L\'evy processes.
Recall that the overall minimum of $(\xi,\pr),$ $-I^{\xi},$ follows
an exponential law of parameter $\gamma>0$ for some $\gamma$ which
is determined in terms of the law $\pr.$ (See Bertoin~\cite{be},
Chapter VII.) Furthermore, it has been proved by
Bertoin~\cite{be:premin} that the pre-minimum part of $(\xi,\pr)$
has the same law as a real valued L\'evy process, say
$(\xi,\pr^{\downarrow}),$ killed at its first hitting time of
$-\mathfrak{e}$ with $\mathfrak{e}$ a r.v. independent of
$(\xi,\pr^{\downarrow})$ and that follows an exponential law of
parameter $\gamma.$ (The process $(\xi,\pr)$ can be viewed as
$(\xi,pr)$ conditioned to drift to $-\infty.$)  The translation of
Bertoin's results for positive self-similar Markov process leads to
the following Proposition. Denote by $\pc,$ the law of the process
obtained by applying Lamperti's transformation to the L\'evy process
$(\xi,\pr^{\downarrow}).$
\begin{prop}\label{prop:1nonegjumps}
If $(X,\p),$ equivalently $(\xi,\pr),$ has no negative jumps, then
there exists a real $\gamma>0$ such that for any $x>0$
$$\p_x(I^X\leq \epsilon)=(\epsilon/x)^{\gamma}\wedge 1,\quad
\epsilon\geq 0,$$ and the law of $\left((X_t,0\leq t<
m),\p_x\right)$ is the same as that of $\left((X_t,0\leq
t<T(Z)),\pcx\right),$ where $Z$ is a random variable independent of
$(X,\pcx)$ and such that $(-\log(Z/x),\pcx)$ follows an exponential
law of parameter $\gamma>0.$
\end{prop}
\begin{proof} This follows from Proposition \ref{pathdecpssmp}
and Theorem~2 in~\cite{be:premin}, described above.
\end{proof}

So to reach our end, we will next provide a description of the
pre-minimum of a real valued L\'evy process that drifts to $\infty,$
which generalizes Bertoin's result and is analogous to the
description of the pre-minimum of a L\'evy process conditioned to
stay positive that has been obtained in \cite{ch2} and
\cite{duquesne2003}.

Let $\widehat{V}(\mathrm{d}x), x\geq 0$ be the renewal measure of
the downward ladder height process, see e.g.~\cite{be} or
~\cite{chaumontdoney} for background. In the remaining of this
Section we will assume that under $\pr,$
\begin{itemize}
\item[(H)]$\begin{cases}& 0 \ \text{is regular for}\ ]-
\infty,0[\\ &\xi\ \text{derives towards}\
+\infty\\ &\text{the measure}\ \widehat{V}(\mathrm{d}x)\ \text{is
absolutely continuous w.r.t Lebesgue's measure.}\end{cases}$
\end{itemize}

In order to construct the L\'evy process which describes the
pre-minimum part of $(\xi,\pr)$ we will need the following Lemma
which is reminiscent of Theorem~1 in \cite{silverstein}. Let
$\pr^{]-\infty,0[}$ be the law of $(\xi,\pr)$ killed at its first
hitting time of $]-\infty,0[.$

\begin{lemma}\label{excdensity}
Under the assumptions $(H)$ the renewal measure
$\widehat{V}(\mathrm{d}x)$ has a density, say $\varphi:\re\to\re^+,$
which is excessive for the semigroup of $(\xi,\pr^{]0,\infty[})$ and
$0<\varphi(x)<\infty$ for a.e. $x\in\re^+.$
\end{lemma}
\begin{proof}
It is known that the processes $(\xi,\pr)$ and $(-\xi,\pr)$ are in
weak duality w.r.t. Lebesgue's measure, so by Hunt's switching
identity we have that $(\xi,\pr^{]-\infty,0[})$ and
$(-\xi,\pr^{]-\infty,0[})$ are also in weak duality w.r.t.
Lebesgue measure, see e.g. \cite{getoorsharpe}. On the other hand,
it is known that the measure $\widehat{V}(dx)$ is an invariant
measure for the process, $S-\xi=\{\sup_{s\leq t}\xi_s -\xi_t,
t\geq 0\},$ $\xi$ reflected at its supremum, see e.g. \cite{be}
Chapter VI exercise 5. So the measure, $\widehat{V}(\mathrm{d}x)$
is excessive for $S-\xi,$ killed at its first hitting time of $0,$
so for $(-\xi,\pr^{]-\infty,0[}).$ Thus the first assertion of
Lemma~\ref{excdensity} is a direct consequence of Theorem
in~Chapter XII paragraph 71 in \cite{dellacheriemeyer}. To prove
the second assertion we recall that
$$\widehat{V}[0,x]=k\pr(-\inf_{0\leq
 s<\infty}\xi_s\leq x), \qquad x\geq 0,$$ with $k\in]0,\infty[$ a
 constant, see \cite{be}
 Proposition~VI.17. So $\varphi<\infty$ a.e. and by the
regularity for $]-\infty,0[$ of 0, the support of the law of
$\inf_{0\leq s<\infty}\xi_s$ is $]-\infty,0[,$ thus $0<\varphi$ a.e.
\end{proof}

Let $\prcp,$ be the  $h$-transform of the law, $\pr^{]-\infty,0[},$
via the excessive function $\varphi.$ That is, $\prcp$ is the unique
measure which is carried by $\{0<\zeta\}$ and under which the
canonical process is Markovian with semi-group $(P^{\searrow}_t,
t\geq 0),$
$$P^{\searrow}_tf(x)=\begin{cases}\frac{1}{\varphi(x)}\er^{]-
\infty,0[}_x\left(f(\xi_t)\varphi(\xi_t)\right)&\ \text{if}\
x\in\{z\in\re: 0<\varphi(z)<\infty\},\\ 0 &\ \text{if}\ x\notin
\{z\in\re: 0<\varphi(z)<\infty\}.\end{cases}$$ Let
$\Lambda=\{z\in\re: 0<\varphi(z)<\infty\}.$ Furthermore, the measure
$\prcp$ is carried by $\{\xi_t \in \Lambda, \xi_{t-}\in\Lambda,
t\in]0,\zeta[\},$ and for any $\mathcal{G}_t$-stopping time $T$
$$\sideset{}{_x^{\searrow}}\pr 1_{\{T<\zeta\}}=
\frac{\varphi(\xi_{T})}{\varphi(x)}1_{\{T<\zeta\}}\pr^{]-\infty,0[}_x,
\qquad \text{on}\ \mathcal{G}_T.$$

In the case where the semigroup of $(\xi,\pr)$ is absolutely
continuous,  $\prcp$ has been introduced in~\cite{ch2} where it is
proved that this measure can be viewed as the law of $(\xi,\pr)$
conditioned to hit 0 continuously. In the case where $(\xi,\pr)$
creeps downward $\varphi$ can be made explicit:
$$\varphi(x)=c\pr(\xi_{T_{]-\infty,-x[}}=-x)>0,\qquad x>0,$$ with
$0<c<\infty,$ a constant, see \cite{be}~Theorem VI.19, and then we
have the right conditioning:
$$\sideset{}{_x^{\searrow}}\pr=\pr^{]-\infty,0[}_x(\ \cdot\ |\,
\xi_{T_{]-\infty,0[}}=0).$$
So in the sequel we will refer to $\prcp$ as the law of $(\xi,\pr)$
conditioned to hit $0.$

\begin{lemma}\label{lemma:descpremin}
Let $\xi$ be a real valued L\'evy process that satisfies the
hypotheses $(H)$ and $\varphi$ be the density of the renewal measure
$\widehat{V}$ as in Lemma~$\ref{excdensity}$. Then for any bounded
measurable functional $F$,
\begin{equation*}\label{eq:condhit0continuosly}
\er\left(F\left(\xi_s-\xi_{\rho-}, 0\leq
s<\rho\right)\right)=\frac{1}{\widehat{V}]0,\infty[}\int_{]0,\infty[}\mathrm{d}a
\varphi(a)\sideset{}{_a^\searrow}\er\left(F\left(\xi_s, 0\leq s<
\zeta\right)\right).
\end{equation*} In particular under
$\pr$ conditionally on $I^{\xi}=a$, the pre-minimum process has the
same law as $\xi+a$ under $\sideset{}{_{-a}^{\searrow}}\pr.$
\end{lemma}
Observe that Bertoin's~\cite{be:premin}~Theorem 2 can be deduced
from this Lemma since in the case where $\xi$ has no negative jumps
$\varphi,$ is given by $\varphi(x)=\gamma e^{-\gamma x}, x>0,$ and
so we have that
\begin{equation*}
\begin{split}
\er\left(F\left(\xi_s, 0\leq s<\rho\right)\right)
&=\int_{]0,\infty[}\mathrm{d}a
 \gamma e^{-\gamma a} \sideset{}{_a^\searrow}\er\left(F\left(\xi_s-a, 0\leq s<
\zeta\right)\right)\\
&=\int_{]0,\infty[}\mathrm{d}a
 \gamma e^{-\gamma a} \er^{]-\infty,0[}_a\left(F\left(\xi_s-a, 0\leq s<
\zeta\right) |\ T_{]-\infty,0[}<\infty\right)\\
&=\gamma \int_{]0,\infty[}\mathrm{d}a \er\left(F\left(\xi_s, 0\leq
s< T_{]-\infty,-a[}\right), T_{]-\infty,-a[}<\infty\right)\\
&=\int_{]0,\infty[}\mathrm{d}a \gamma e^{-\gamma a}
\er\left(F\left(\xi_s, 0\leq s< T_{]-\infty,-a[}\right)
e^{-\gamma \xi_{T_{]-\infty,-a[}}},\ T_{]-\infty,-a[}<\infty\right)\\
&=\erc\left(F\left(\xi_s, 0\leq s<
T_{]-\infty,-\mathfrak{e}[}\right)\right),
\end{split}
\end{equation*}
where $\prc$ and $\mathfrak{e}$ are as explained just before
Proposition~\ref{prop:1nonegjumps}.
\begin{proof}
To prove the claimed identity, we will start by calculating for
any continuous and bounded functional $F,$
$$\er^{\mathrm{e/\lambda}}\left(F\left( \xi_s-\xi_\rho, 0\leq
s<\rho\right)\right),$$ where $\er^{\mathrm{e/\lambda}}$ is the law
of $(\xi,\pr)$ killed at time $\mathrm{e}/\lambda$, with
$\mathrm{e}$ an exponential random variable independent of
$(\xi,\pr).$ To do that we will denote by $\{\underline{L}_u, u\geq
0\}$ the local time at 0 of the strong Markov process $\{\xi_t-I_t,
t\geq 0\},$ by $g_t$ the last hitting time of 0 by $\xi-I$ before
time $t,$ $g_t=\sup\{s\leq t: \xi_s-I_s=0\},$ and by $\underline{N}$
the excursion measure of $\xi-I$ away from 0. Indeed, using
Maisonneuve's exit formula of excursion theory it is justified that
\begin{equation*}
\begin{split}
\er^{\mathrm{e/\lambda}}\left(F\left( \xi_s-I_{\rho}, 0\leq
s<\rho\right)\right)& = \int^{\infty}_0\mathrm{d}t\lambda
e^{-\lambda t}\er\left(F\left( \xi_s-I_{g_t}, 0\leq s<g_t
\right)\right)\\
&=\int^{\infty}_0\mathrm{d}t\lambda e^{-\lambda
t}\er\left(\int^t_0\mathrm{d}\underline{L}_u F\left( \xi_s-I_{u},
0\leq s<u\right)\underline{N}(t-u<\zeta)\right)\\
&=\er\left(\int^\infty_0\mathrm{d}\underline{L}_u e^{-\lambda u}
F\left(\xi_s-I_{u}, 0\leq
s<u\right)\right)\underline{N}(1-e^{-\lambda \zeta}).
\end{split}
\end{equation*}
Next, making $\lambda$ tend to 0, the left hand term in the previous
equality tends to $$\er\left(F\left( \xi_s-\xi_{\rho-}, 0\leq
s<\rho\right)\right),$$ while the right hand term tends to
$$\er\left(\int^\infty_0\mathrm{d}\underline{L}_u F\left(\xi_s-I_{u}, 0\leq
s<u\right)\right)\underline{N}(\zeta=\infty).$$ Finally, a
straightforward extension of Lemma~3 in~\cite{ch3} to our weaker
hypothesis allows us to ensure that
$$\er\left(\int^\infty_0\mathrm{d}\underline{L}_u F\left(\xi_s-I_{u}, 0\leq
s<u\right)\right)=\int_{]0,\infty[}\mathrm{d}a
\varphi(a)\sideset{}{_a^\searrow}\er\left(F\left(\xi_s, 0\leq
s<\zeta\right)\right),$$ which concludes the proof given that
$\widehat{V}]0,\infty[=\left(\underline{N}(\zeta=\infty)\right)^{-1}.$
\end{proof}

We next introduce the law of a L\'evy process conditioned to hit by
above a given level $a\in\re$. Owing to the fact that the function
$\varphi$ is excessive for $(\xi,\pr)$ killed at 0, we have that for
any $a\in\re$ the function $\varphi_a:\re\to\re^+$ defined by
$\varphi_a(x)=\varphi(x-a), x\in\re,$ is excessive for the semigroup
of $\xi$ killed at its first hitting time of $]-\infty,a[.$ Indeed,
\begin{equation*}
\er_x(\varphi_a(\xi_t),t<T_{]-\infty,a[})=
\er_{x-a}(\varphi(\xi_t),t<T_{]-\infty,0[})\leq
\varphi(x-a)=\varphi_a(x),\qquad x>a,
\end{equation*}
and analogously it is verified that $\lim_{t\to 0+}
\er_x(\varphi_a(\xi_t),t<T_{]-\infty,a[})=\varphi_a(x).$ We will
denote by $\sideset{}{^{\searrow a}}\pr$ the $h$-transform of the
law of $\xi$ killed at it first hitting time of $]-\infty,a[$ via
$\varphi_a$, i.e.: for any $\mathcal{G}_t$-stopping time $T$, with
an obvious notation for $\pr^{]-\infty,a[}_x$,
$$\sideset{}{_x^{\searrow a}}\pr 1_{\{T<\zeta\}}=
\frac{\varphi_a(\xi_{T})}{\varphi_a(x)}1_{\{T<\zeta\}}\pr^{]-\infty,a[}_x,
\qquad \text{on}\ \mathcal{G}_T.$$ The following elementary Lemma
will enable us to refer to this measure as the law of $\xi$
conditioned to hit $a$ continuously and by above. Of course the
measure $\sideset{}{_x^{\searrow 0}}\pr$ is simply
$\sideset{}{_x^{\searrow}}\pr.$

\begin{lemma}\label{lemma:condhitcontinuously}
Let $(\xi,\pr)$ be a real valued L\'evy process that satisfies the
hypotheses $(H)$. For $a\in\re,$ and any $x>a$ the law of $\xi+a$
under $\sideset{}{_{x-a}^{\searrow}}\pr$ is the same as that of
$\xi$ under $\sideset{}{_x^{\searrow a}}\pr.$ As a consequence, for
a.e. $x>a$
$$\sideset{}{_x^{\searrow a}}\pr\left(\xi_0=x;\ \zeta<\infty;\
\xi_t>a \ \text{for all}\ t<\zeta;\ \xi_{\zeta-}=a \right)=1.$$
\end{lemma}
\begin{proof}To prove the first assertion it suffices to
verify that both laws are equal over $\mathcal{G}_t$ for any $t>0.$
Indeed, the spatial homogeneity of $(\xi,\pr)$ implies that for
$t>0$  and any bounded measurable functional $F$
\begin{equation*}
\begin{split}
\sideset{}{_x^{\searrow a}}\pr\left(F(\xi_s,0\leq s<t)1_{\{t<
\zeta\}}\right)
&=\frac{1}{\varphi_a(x)}\pr_x(F(\xi_s,0\leq s<t)
1_{\{t<T_{]-\infty,a[}\}}\varphi_a(\xi_t))\\
&=\frac{1}{\varphi(x-a)}\pr_{x-a}(F(\xi_s+a,0\leq s<t)
1_{\{t<T_{]-\infty,0[}\}} \varphi(\xi_t))\\
&=\sideset{}{_{x-a}^{\searrow}}\pr(F(\xi_s+a,0\leq
s<t)1_{\{t<\zeta\}}).
\end{split}
\end{equation*}
Now, the second assertion is an easy consequence of
Lemma~\ref{lemma:descpremin} and the hypothesis that $0$ is regular
for $]-\infty,0[.$
\end{proof}
A rewording of Lemma~\ref{lemma:descpremin} using
Lemma~\ref{lemma:condhitcontinuously} reads:

\begin{teo}\label{preminncond}
Let $(\xi,\pr)$ be a real valued L\'evy process that satisfies the
hypotheses $(H)$. The following identity holds for any bounded
measurable functional $F$,
\begin{equation}\label{eq:descpreminimumgen}
\er_{x}\left(F\left(\xi_s, 0\leq
s<\rho\right)\right)=\frac{1}{\widehat{V}]0,\infty[}\int_{]-\infty,x[}\mathrm{d}a
\varphi_a(x)\sideset{}{_x^{\searrow a}}\er\left(F\left(\xi_s, 0\leq
s< \zeta\right)\right).
\end{equation}
\end{teo}
We have now all the elements to state the main result of this
section whose proof follows easily from
Lemma~\ref{lemma:condhitcontinuously} \& Theorem~\ref{preminncond}.

\begin{teo}
Let $(\xi,\pr)$ be a real valued L\'evy process that satisfies the
hypotheses $(H)$ and $(X,\p)$ be the self-similar Markov process
associated to $(\xi,\pr)$ via Lamperti's representation. Then for
any bounded measurable functional $F$,
\begin{equation*}\begin{split}
\e_x\left(F\left(X_s,0\leq s<m\right)\right)
&=\int^1_0\nu_1(\mathrm{d}v)\sideset{}{_x^{\searrow
v x}}\p(F(X_s, 0\leq s<\zeta))\\
&=\int^1_0\nu_1(\mathrm{d}v)\sideset{}{_{1/v}^{\searrow 1}}\p(F(v x
X_{s/vx}, 0\leq s<vx\zeta)),
\end{split}
\end{equation*}
where $\nu_1$ is a measure over $]0,1[$ with density
$$\frac{\nu_1(\mathrm{d}v)}{\mathrm{d}v}=
(\widehat{V}]0,\infty[)^{-1}v^{-1}\varphi(-\ln v),\qquad 0<v<1,$$
and $\sideset{}{_x^{\searrow vx}}\p$ is the law of the process
obtained by applying Lamperti's representation to
$(\xi,\sideset{}{_{\log(x)}^{\searrow \log(vx)}}\pr).$
\end{teo}
\end{section}
\begin{section}{The asymptotic behavior of the pre- and post-minimum
as the minimum tends to 0.}\label{prepostas} Throughout this section
we will leave aside the assumptions $(H)$. We only assume that the
underlying L\'evy process $\xi$ drifts to $\infty$, it is not a
subordinator and it is non lattice. Some ancillary hypothesis will
be stated below.

\begin{subsection}{Post-minimum}

Under these hypotheses, it is known that the support of the law of
$I^\xi$ is $]-\infty,0]$. From (3), the support of $I^X$ is then
$[0,y]$ under $\p_y$, for any $y>0$. Proposition~\ref{pathdecpssmp}
shows that a regular version of the law of the post-minimum process
$(X_{m+t},\,t\ge0)$ under $\p_y$ given $I^{X}=x$, for $x\in]0,y]$ is
given by the law of the process $\left(\left(x
\poxi_{\;\taup\,(t/x)},\,t\ge0\right),\pr\right)$. In particular,
this law does not depend on $y$. Let us denote it by $\pp^x$. A
straight consequence of this representation is that the family
$(\pp^x)$ is weakly continuous on $]0,\infty[$. In
Theorem~\ref{teo:postminxto0} below, we show that if moreover
$0<\er(\xi_1)<\infty,$ then $\pp^x$ converges weakly as $x$ tends to
0 towards the law $\p_{0+}$.  This measure is the weak limit of
$\p_x$ as $x\to 0+$, whose existence is ensured by Theorem~2 in
\cite{chaumontcaballero}.

Recall that Millar's results implies that for any $x>0,$ the process
$(X,\pp^x)$ is strongly Markov with values in $[x,\infty[$.\\

\begin{teo}\label{teo:postminxto0}
Assume that $0<\er(\xi_1)<\infty.$ The laws $\ppx$ converge weakly
in $\mathbb{D}$ as $x\to 0+$ to the law $\p_{0+}.$ As a consequence,
for any $x>0,$ $$\p_x(\cdot \circ\theta_{m}|
I^{X}<\epsilon)\xrightarrow[\epsilon\to 0]{w}\p_{0+}(\cdot).$$
\end{teo}
\begin{proof}
Recall that from~\cite{chaumontcaballero}, under our hypothesis, the
family of laws $(\p_x)$ converges weakly in $\mathbb{D}$ as
$x\downarrow0$ towards the non degenerate law of a self-similar
strong Markov process. Denote by $\p_{0+}$ the limit law. Then on
the space $\mathbb{D}$, we define a process $X^{(0)}$ with law
$\p_{0+}$. We recall from \cite{chaumontcaballero} that
\begin{equation}\label{limites}
\lim_{t\rightarrow 0+}X^{(0)}_t=0
\;\;\;\mbox{and}\;\;\;\lim_{t\uparrow\infty}X_t^{(0)}=+\infty
\,,\;\;\; \p_{0+}\ \mbox{a.s.}
\end{equation}
Let $(x_n)$ be any sequence of positive real numbers which tends
to 0. Define $\Sigma_n=\inf\{t:X^{(0)}_t\ge x_n\}$, then by the
Markov property and Lamperti's representation, we have
\begin{equation}\label{Lamprandom}
Y^{(n)}\eqdef(X^{(0)}_{\Sigma_n+t},\,t\ge0)=
\left(X_{\Sigma_n}^{(0)}\exp\xi^{(n)}_{\tau^{(n)}(t/X_{\Sigma_n}^{(0)})},\,t\ge0\right)\,,
\end{equation}
where on the left hand side of the second equality,
$X_{\Sigma_n}^{(0)}$ and $\xi^{(n)}$ are independent and
$\xi^{(n)}\ed\xi$. Let
\[I_n=\inf_{t\ge0}Y^{(n)}_t\;\;\;\mbox{and}\;\;\;m_n=\sup\{t:Y_t^{(n)}\wedge
Y_{t-}^{(n)}=I_n\}\,.\] Then we deduce from (\ref{Lamprandom}) and
Proposition 1 the following representation:
\begin{equation}\label{postmin1}
(Y^{(n)}_{m_n+t},\,t\ge0)=
\left(I_n\exp\poxi^{(n)}_{\taup(t/I_n)},\,t\ge0\right)\,,
\end{equation} where $\poxi^{(n)}$ is independent of the events prior to
$m_n$. In particular, $\poxi^{(n)}$ is independent of ${\cal
G}_n\eqdef\sigma\{I_k:k\ge n\}$. It follows from (\ref{postmin1})
that for any bounded and measurable functional $H$,
\begin{equation}\label{postmin2}
\e_{0+}(H(Y^{(n)}_{m_n+t},\,t\ge0)\,|\,{\cal G}_{n})=\ep^{I_n}(H)\,.
\end{equation}
Since $(X,\p_x)$, $x\ge0$ is a Feller process, the tail
$\sigma$-field $\cap_{t>0}\sigma\{X^{(0)}_s:s\le t\}$ is trivial and
it is not difficult to check that $\cap_n {\cal G}_n\subset
\sigma\{X^{(0)}_s:s\le t\}$ for each fixed $t$. So $\cap_n {\cal
G}_n$ is trivial. On the other hand, from (\ref{limites}) we have
$\lim_n\Sigma_n=0$ and $\lim_n m_n=0,$ $\p_{0+}$--a.s., so
\[(Y^{(n)}_{m_n+t},\,t\ge0)\longrightarrow X^{(0)}\,,\;\;\; \p_{0+}\ \mbox{a.s., as
$n\rightarrow+\infty$,}\] on the space $\mathbb{D}$. Hence if we
suppose moreover that $H$ is continuous, then
\begin{equation}\label{convps1}
\lim_nE(H(Y^{(n)}_{m_n+t},\,t\ge0)\,|\,{\cal
G}_{n})=\lim_n\ep^{I_n}(H)=\e_0(H)\,,\;\;\; \p_{0+}\ \mbox{almost
surely}\,.
\end{equation}
Now, from (\ref{time}), we have $I_n=X_{\Sigma_n}^{(0)}\exp
I^{\xi^{(n)}}$. Recall that from Theorem 1 in [11], the r.v.
$X_{\Sigma_n}^{(0)}$ may be decomposed as $X_{\Sigma_n}^{(0)}\ed
x_ne^{\theta}$, where $\theta$ is a finite r.v. whose law is this of
the limit overshoot of the L\'evy process $\xi$, i.e. if
$T_z=\inf\{t:\xi_t\ge z\}$, then under our hypothesis, $\xi_{T_z}-z$
converges in law as $z\uparrow+\infty$ towards the law of $\theta$.
So, we have
\begin{equation}\label{idloi}
I_n\ed x_n e^\theta e^{I^\xi}\,,\end{equation} where $\theta$ and
$I^\xi$ are independent. On the space ${\mathbb{D}}$, we define a
r.v. $\nu$ such that $\nu\ed e^\theta e^{I^\xi}$ (so that $I_n\ed
x_n\nu$), then it follows from (\ref{convps1}) that
\begin{equation}\label{convprob}
\ep^{x_n\nu}(H)\longrightarrow\e_0(H)\,,\;\;\;\mbox{in
probability, as $n\rightarrow+\infty$}\,.
\end{equation}
So there exists a subsequence $x_{n_k}$ such that
\begin{equation}\label{convps2}
\ep^{x_{n_k}\nu}(H)\longrightarrow\e_0(H)\,,\;\;\;\mbox{a.s., as
$k\rightarrow+\infty$}\,.
\end{equation}
The convergence (\ref{convps2}) implies that there exists
$\omega_0\in{\mathbb{D}}$ such that $\nu(\omega_0)>0$ and
$\ep^{x_{n_k}\nu(\omega_0)}(H)\rightarrow\e_0(H)$, as
$k\rightarrow+\infty$. Put $a=\nu(\omega_0)$ and for all
$\omega\in{\mathbb{D}}$ define
$S_a(\omega)=(a^{-1}\omega_{at},\,t\ge0)$. Since $S_a$ is a
continuous functional on ${\mathbb{D}}$, we have
\[
\ep^{x_{n_k}a}(H\circ S_a)\longrightarrow\e_0(H\circ
S_a)\,,\;\;\;\mbox{as $k\rightarrow+\infty$}\,. \] But from the
scaling property, we have $\ep^{x_{n_k}a}(H\circ
S_a)=\ep^{x_{n_k}}(H)$ and $\e_{0+}(H\circ S_a)=\e_{0+}(H)$. In
conclusion, for any bounded and continuous functional $H$ on
${\mathbb{D}}$ and for any sequence $(x_n)$ which decreases to 0,
there is a subsequence $(x_{n_k})$ such that
$\ep^{x_{n_k}}(H)\longrightarrow\e_{0+}(H)$, as $k$ tends to
$\infty$. This proves our result.
\end{proof}
\end{subsection}
\begin{subsection}{Pre-minimum}
In our description of the pre-minimum process we have provided,
under some assumptions, a method to construct a process that can be
viewed as $X$ conditioned to die at a given level $0<a<X_0.$ But a
priori this method cannot be applied to construct a process that
dies at 0, since this means conditioning the underlying L\'evy
process to die at $-\infty.$ Thus, the purpose of this section is to
construct the law of a self-similar Markov process that can be
viewed as the law of a pssMp that drift to $\infty,$ $X,$
conditioned to hit 0 in a finite time. In fact we will answer the
questions: \emph{What is the process obtained by making tend to $0$
the value of the overall minimum of $X$? Is the resulting process
determined by the pre-minimum process of $X$?} In the case where $X$
 has no negative jumps, using the assertions in
Proposition~\ref{prop:1nonegjumps} it is clear, at least
intuitively, that the process $(X,\pc)$ can be obtained from
$(X,\p)$ by making tend to 0 the value of its overall minimum.
Actually, the former process can be viewed as $(X,\p)$ conditioned
to have an overall minimum equal to 0 and this conditioning depends
only on the pre-minimum part of $(X,\p).$

As a consequence of the assumption that $(X,\p)$ drifts to $\infty$,
the set of paths that have an overall minimum equal to 0,
$\{I^{X}=0\},$ has probability 0, and so the law of $(X,\p)$
conditionally on that set does not make sense. A natural issue to
give a meaning to that conditioning is by approximating that set by
the sequence $\{I^{X}<\epsilon\}$ as $\epsilon \to 0.$ So, our main
task will be describe the limit law of the pre-minimum process
conditionally on the event $\{I^{X}<\epsilon\}$ as $\epsilon\to 0$.
To that end we will use the method of $h$-transformations.

Let $h:]0,\infty[\to[0,\infty]$ be the function defined by
\begin{equation}\label{eq:ct0} h(x)=\liminf_{\epsilon\to 0}
\frac{\p_x(I^X< \epsilon)}{\p_1(I^X< \epsilon)},\qquad x\in ]0,\infty[.
\end{equation}
The following Lemma will be useful.
\begin{lemma}
The function $h$ defined in equation~$(\ref{eq:ct0})$ is excessive
for the semigroup of the pssMp $X.$
\end{lemma}

\begin{proof}Given that the cone of excessive functions is closed under
$liminf$ it suffices with proving that for every $\epsilon>0,$
the function
\begin{equation*} h^\epsilon(x)=\frac{\p_x(I^X< \epsilon)}
{\p_1(I^X< \epsilon)},\qquad x\in ]0,\infty[,
\end{equation*}
is excessive for the semigroup of $X.$ Indeed, owing the relation
$$\p_x(I^X< \epsilon)=\p_x(L_{\epsilon}>0),\qquad \text{with}\
L_{\epsilon}=\sup\{s>0 : X_s < \epsilon\},\ (\sup\{\emptyset\}=0),$$
and the Markov property, it is straightforward that for any reals
$t>0$ and $x>0$
\begin{equation*}
P_th^\epsilon(x)=\frac{\e_x(\p_{X_t}(L_{\epsilon}>0))}{\p_1(I^X<
\epsilon)}=\frac{\p_x((L_{\epsilon}-t)^+>0)}{\p_1(I^X<\epsilon)}\leq
h^\epsilon(x),\qquad x>0,
\end{equation*}
and  $$\lim_{t\to0}P_th^\epsilon(x)=h^\epsilon(x),\qquad
x>0.$$\end{proof}

To perform the desired conditioning we will make some assumptions
on the excessive function $h$. Firstly, to avoid pathological
cases we will assume that $h$ does not take the values 0 or
$\infty,$ and  next that it has some regularity, namely that
\begin{enumerate}\item[\bf{(H')}] the liminf
in equation~(\ref{eq:ct0}) is in fact a limit and
$h:]0,\infty[\to]0,\infty[$ is a non-constant function.
\end{enumerate}
The hypothesis (H') is satisfied by a wide class of positive
self-similar Markov processes, as it will be seen in
Remark~\ref{example1} below, and, whenever it holds, the
self-similarity implies that, the excessive function $h$ has the
form
$$h(x)=x^{-\gamma},\qquad x>0,\ \text{for some}\ \gamma>0.$$
Here is a reformulation of (H') in terms of the underlying L\'evy
process $(\xi,\pr)$. First, one has $\pr(-I^{\xi} > z)>0$ for each
$z>0$ and
\begin{equation*}
\lim_{u\to\infty}\frac{\pr(-I^{\xi}>u-z)}{\pr(-I^{\xi}>u)}=e^{\gamma
z}\quad\text{for each} \ z\in\re.
\end{equation*}
In other words, the law of the negative of the overall minimum of
$\xi$ belongs to one of the classes $\mathcal{L}^{\gamma},$ for some
$\gamma>0;$

In the sequel we will assume that the hypothesis (H') is satisfied.
Let $\pc$ be the $h$-transform measure of $\p$ via $h$, i.e.: for
any $\mathcal{F}_t$-stopping time $T$
$$\pcx 1_{\{T<\zeta\}}=
\frac{h(X_{T})}{h(x)}\p_x, \qquad \text{on}\ \mathcal{F}_T.$$

By standard arguments it follows that the law $\pc$ is that of a
positive self-similar Markov process, say $(X,\pc).$ We will denote
by $(\xi,\prc)$ the L\'evy process associated to $(X,\pc)$ via
Lamperti's transformation. By the absolute continuity relation
between $\pc$ and $\p$ applied to the sequence of
$\mathcal{F}$-stopping times
$$T_t=\inf\{r>0: \int^r_0X^{-1}_s\mathrm{ds}>t\},\qquad t\geq 0,$$
and Lamperti's transformation, it holds that
$\er(e^{-\gamma\xi_t})\leq 1$ for all $t>0,$ and more importantly
that the laws $\prc$ and $\pr$ are absolutely continuous: for any
$t\geq 0$
\begin{equation}\label{eq:htransflevy}\prc1_{\{t<\zeta\}}=
e^{-\gamma\xi_t}\pr,\quad
\text{on} \ \mathcal{F}_{T_t}=\mathcal{G}_t.\end{equation}  The
latter relation can be extended to $\mathcal{G}$-stopping times
using standard arguments.

With the following result we prove that the family of laws $(\pcx,
x>0)$ can be thought as those of the process $(X, \p)$ strictly
before $m$ when the whole trajectory is conditioned to have an
overall minimum equal to 0.

\begin{teo}\label{convmin0}Assume the hypothesis $(H')$ is satisfied.
\begin{enumerate}
\item[$(i)$]The process $(X,\pcx)$ hits $0$ in a finite time, a.s. Moreover,
$$\pcx\left(T_0<\infty, X_{T_0-}=0\right)=1, \quad \text{for all}\ x>0,$$
if and only if Cram\'er's condition, $\e(e^{-\gamma \xi_1})=1,$ is
satisfied. \item[$(ii)$] If $(\xi,\pr)$ satisfies furthermore that
either
\begin{enumerate}
\item[$(a1)$] its law is not lattice,
\item[$(a2)$] Cram\'er's condition, $\e(e^{-\gamma\xi_1})=1$ and
$\e(\xi^{-}_1e^{-\gamma\xi_1})<\infty$ are satisfied,
\end{enumerate} or
\begin{enumerate}
\item[$(b1)$] $\e(e^{-\gamma\xi_1})<1,$
\end{enumerate} then the law $\pc$ is determined by the law of the
pre--minimum process of $(X,\p)$
in the following way: for any $x>0$
\begin{equation*}
\lim_{\epsilon\to 0+}\p_x(F_t\cap\{t<m\} | I^X < \epsilon)=
\p^{\downarrow}_x(F_t\cap\{t<T_0\}),\quad F_t\in\mathcal{F}_t,\;\;
t\geq 0.
\end{equation*}
\end{enumerate}
\end{teo}
\noindent A consequence of (ii) in Theorem~\ref{convmin0} is that
the finite dimensional laws of the pre-minimum process converge to
those of
$(X,\pcx).$\\

\noindent {\it Proof of part $(i)$}. By the
identity~(\ref{eq:htransflevy}) it follows that
$$\er(e^{-\gamma\xi_t})=\prc(t<\zeta),\qquad\text{for all}\ t>0,$$
and so under $\prc$ the canonical process $\xi$ has an infinite
lifetime if and only if $\er(e^{-\gamma\xi_t})=1,$ for all $t>0$ or
equivalently for some $t>0,$ see e.g. Sato~\cite{sa} Theorem~25.17.
In which case Cram\'er's condition is satisfied and the process
$(\xi,\prc)$ drifts to $-\infty.$ Given that the process $(X,\pc)$
coincides with the pssMp associated to $(\xi,\prc)$ via Lamperti's
transformation, we conclude using Lamperti's representation of
pssMp, see Section~\ref{lamptrasform}, that if Cram\'er's condition
is satisfied then $$\pcx(T_0<\infty, X_{T_0-}=0)=1, \quad \text{for
all}\ x>0.$$ Now, assume that Cram\'er's condition is not satisfied,
that is $\er(e^{-\gamma\xi_t})<1$ for some $t>0.$ By Theorem~25.17
in \cite{sa} this implies that the latter holds for all $t>0.$ So
the L\'evy process $(\xi,\pr)$ has a finite lifetime, actually it is
a real valued L\'evy process that has been killed at an independent
time that follows an exponential law of parameter
$\kappa=-\log\er(e^{-\gamma\xi_1}).$ According to Lamperti
representation of pssMp we have that in this case
$$\pcx(T_0<\infty, X_{T_0-}>0)=1, \quad \text{for
all}\ x>0.$$ In any case, $(X,\pcx)$ hits 0 in a finite time a.s.
for all $x>0$. Which finish the proof of assertion (i).\\

\noindent {\it Proof of part $(ii)$}. To prove the assertion we will
start by proving that for any $x>0$ and $t>0,$
\begin{equation}\label{eq1:convmin0} \lim_{\epsilon\to 0+}\p_x(t < m
| I^X<\epsilon)=
x^{\gamma}\e_x(X^{-\gamma}_t)=\sideset{}{_x^\downarrow}\p(t<T_0).
\end{equation}
To that end we will use that $\{L_{\epsilon}>0\}=\{I^X <
\epsilon\},$ and so that
\begin{equation*}
\begin{split}
\p_x\left(t < m, I^X<\epsilon\right)&=
\p_x\left(t<m, 0<L_{\epsilon},t<L_{\epsilon}\right)\\
&=\p_x\left(m \wedge L_{\epsilon} > t\right)\\
&=\p_x\left(\p_{X_t}\left(m \wedge L_{\epsilon} >
0\right)\right)\\
&=\p_x\left(\p_{X_t}\left(L_{\epsilon} > 0\right)\right),
\end{split}
\end{equation*}
which is a consequence of the fact that $L_{\epsilon}$ and $m$ are
both coterminal times, the Markov property and that $\p_x(m=0)=0,$
owing that $(\xi,\pr)$ is not a subordinator. Moreover, it follows
from the scaling and Markov properties that
$$\e_x\left(\p_{X_t}\left(L_{\epsilon} > 0\right)\right)=
\e_x\left(g(X_t/\epsilon)\right),$$ where $g(z)=\p_{x}\left(I^X \leq
z^{-1}\right).$ Now, if the conditions (a-1,2) are satisfied, then
the main result of \cite{bd} implies that $g(z)=z^{\gamma}L(z),
z>0,$ where $L:]0,\infty[\to]0,\infty[$ is a bounded and slowly
varying function such that $L(z)\xrightarrow[]{}C\in]0,\infty[$ as
$z\to\infty.$ In this case, the dominated convergence theorem
implies that
\begin{equation*}
\begin{split}\lim_{\epsilon\to 0}\p_x(t<m | I^{X} < \epsilon)
&=\lim_{\epsilon\to0}\frac{1}{\p_x(I^X<\epsilon)}\e_x(g(X_t/\epsilon))\\
&=\lim_{\epsilon\to 0}\left(\frac{\epsilon^{\gamma}} {\p_{x}(I^X <
\epsilon)}\right)\e_x\left(X^{-\gamma}_t
L(\epsilon/X_t)\right)\\&=x^{\gamma}\e_x(X^{-\gamma}_t).
\end{split}
\end{equation*}
However, in the case where Cram\'er's condition is not satisfied it
follows from hypothesis (H') that $g$ is regularly varying at
infinity with index $\gamma$ and we claim that
$\e_x(X^{-\gamma-1}_s)<\infty$ for $x>0, t\geq 0,$ which, in view of
Proposition 3 in \cite{breiman}, imply that
$$\lim_{\epsilon\to0}\frac{1}{g(1/\epsilon)}\e_x(g(X_t/\epsilon))=
\e_x(X^{-\gamma}_t),$$
and the limit in equation~(\ref{eq1:convmin0}) follows. So we just
have to prove that $\e_x(X^{-\gamma-1}_s)<\infty$ for $x>0, t\geq
0.$ Indeed, we have seen that hypothesis (H') implies that
$\er(e^{\gamma\xi_t})\leq 1,$ for all $t\geq 0,$ and since
Cram\'er's condition is not satisfied the latter inequality is a
strictly one. So, by Lamperti's transformation
\begin{equation*}
\begin{split}
\int^{\infty}_0\mathrm{d}t\e_x(X^{-\gamma-1}_t)&=
x^{-(\gamma+1)}\er\left(\int^{\infty}_0\mathrm{d}t
\exp\{-(\gamma+1)\xi_{\tau(tx^{-1})}\}\right)\\
&=x^{-\gamma}\int^{\infty}_0\mathrm{d}s\er(e^{-\gamma\xi_s})\\
&=x^{-\gamma}(-\log(\er(e^{-\gamma\xi_1})))<\infty,\qquad x>0.
\end{split}
\end{equation*}
Thus for $x>0,$ $\e_x(X^{-\gamma-1}_t)<\infty,$ for a.e. $t>0,$ and
by the scaling property the latter holds for any $t>0,x>0.$

To conclude, let $F_t\in\mathcal{F}_t,$ $t>0$, then arguing as
before and using Fatou's lemma we have that
\begin{equation*}
\begin{split}
\liminf_{\epsilon \to 0}\p_x(F_t\cap \{t<m\}|I^X < \epsilon)&=
\liminf_{\epsilon \to 0}\left(\frac{\p_1(I^X < \epsilon)}
{\p_x(I^X < \epsilon)}\right)\e_x\left(1_{F_t}\frac{\p_{X_t}
(I^X < \epsilon)}{\p_{1}(I^X < \epsilon)}\right)\\
& \geq x^{\gamma}\e_x(1_{F_t} X^{-\gamma}_t).
\end{split}
\end{equation*}
Furthermore, applying this estimate to the set complementary of
$F_t$ and using the result in equation~(\ref{eq1:convmin0})
we get that
\begin{equation*}
\begin{split}
\limsup_{\epsilon \to 0}\p_x(F_t\cap \{t<m\}|I^X < \epsilon)& \leq
x^{\gamma}\p_x(F_t X^{-\gamma}_t).
\end{split}
\end{equation*}\qed

It is interesting to note that in the non-Cram\'er case the law
$\pc$ is that of a pssMp that hits 0 in finite time and it does it
by a jump,
$$\pcx(T_0<\infty , X_{T_0-}>0)=1,\qquad \forall x>0.$$
Roughly speaking, Theorem~\ref{convmin0} tells us that in this case
by pulling down the trajectory of $(X,\p),$ under the law
$\p_{\cdot},$ from the place at which it attains its overall infimum
for the last time, we break this trajectory and introduce a jump to
the level 0.

However, the equality in (ii) Theorem~\ref{convmin0} does not hold
on the whole $\sigma$-field of the events prior to $m$, i.e.
$\mathcal{F}_{m-}=\sigma\left(F_t\cap\{t<m\}, F_t\in\mathcal{F}_t,
t\geq 0\right)$. Indeed, if this were the case it would imply that
$$\lim_{\epsilon\to 0}\p_x(X_{m-}\in\mathrm{d}y| I^X<\epsilon)
=\p^{\downarrow}_x(X_{T_0-}\in\mathrm{d}y),$$ given that $X_{m-}$ is
$\mathcal{F}_{m-}$--measurable. But the r.h.s. in the previous
equality is equal to $\prc(x \exp\{\xi_{\mathrm{e}}\}\in
\mathrm{d}z),$ where $\mathrm{e}$ is a random variable independent
of $\xi^{\downarrow}$ and with an exponential law of parameter
$\kappa=-\log(e^{-\gamma\xi_1}).$ While the l.h.s. is equal to the
Dirac mass at $0$ whenever $0$ is regular for $(-\infty,0).$

\begin{remark}\label{example1}
Owing to the equivalent formulation of hypothesis (H') in terms of
the underlying L\'evy process it is easy to provide examples of
pssMp that satisfies (H'). Indeed, it is easily deduced from
Proposition~\ref{prop:1nonegjumps} that when the process has no
negative jumps the function $\varphi$ has the properties required in
(H'). Besides, if a L\'evy process does satisfies the hypotheses
(a1) and Cram\'er's condition in (a2) of Theorem~\ref{convmin0}, it
follows from the result of Bertoin and Doney~\cite{bd} that
$$\lim_{t\to\infty} e^{\gamma t}\pr(I^{\xi}<-t)=C,$$ where $C<\infty$
and $C>0$ if and only if $\er(\xi^-_1e^{-\gamma \xi_1})<\infty.$
We deduce therefrom that under (a1) and (a2) of
Theorem~\ref{convmin0} we have
$$\p_x(I^X<\epsilon) \sim \epsilon^{\gamma}x^{-\gamma}C,\qquad
\text{as}\ \epsilon \to 0,$$ and hence (H') is satisfied.
Furthermore, the hypothesis (H') holds if the distribution of the
negative of the overall minimum of $(\xi,\pr)$ belongs to a class of
close to exponential laws $\mathcal{S}^\gamma$ with $\gamma>0$. (See
the recent work~\cite{KKM} for the definition of the classes
$\mathcal{S}^\gamma$ and NASC on the L\'evy process $(\xi,\pr)$ that
ensure that the negative of the overall infimum belongs to one of
this classes.)
\end{remark}
\end{subsection}
\end{section}

\begin{section}{Conditioning a pssMp to hit 0
continuously}\label{condpssMphito0cont} Throughout this section we
will assume that $(X,\p)$ is a self--similar Markov process that
belongs to the class (LC1). It was showed by Lamperti~\cite{la} that
under these assumptions the process $(X,\p)$ is the exponential of a
L\'evy process that has been killed at an independent exponential
time and time changed, see Section~\ref{lamptrasform} for more
details. So, for notational convenience we will hereafter assume
that $(\xi,\pr)$ is a L\'evy process (with infinite lifetime), that
$\boldsymbol{e},$ is an independent r.v. that follows an exponential
law of rate $\ke>0,$ and that the L\'evy process with finite
lifetime associated to $(X,\p)$ via Lamperti's transformation is the
one obtained by killing $(\xi,\pr)$ at time $\boldsymbol{e}.$

The problem of conditioning a self--similar Markov process that hits
0 by a jump to hit 0 continuously is a problem that was studied by
Chaumont~\cite{ch2} in the case where the process has furthermore
stationary and independent increments, i.e. is a stable L\'evy
process. See Chaumont and Caballero~\cite{chaumontcaballero2005} for
a computation of the underlying L\'evy process of this pssMp in
Lamperti's representation.

Throughout this section we  will assume that
\begin{itemize}
\item[(H")]$=\begin{cases} \text{non--arithmetic}\\
\text{there exists a}\ \gamma<0\
\text{for which}\ \er(e^{\gamma\xi_1})=e^{\ke},\\
\er(\xi^-_1e^{\gamma\xi_1})<\infty.\end{cases}$
\end{itemize}
Under these hypotheses we will prove the existence of a
self--similar Markov process $(X,\pc)$ that can be thought as
$(X,\p)$ conditioned to hit 0 continuously.

The second hypothesis in (H") implies that the function
$h^{\scriptscriptstyle{\downarrow}}(x)=e^{\gamma x},x\in\re$ is an
invariant function for the semigroup of $(\xi,\pr),$ killed at time
$\boldsymbol{e}.$ Let $\prc$ be the $h$--transform of the law of
$(\xi,\pr)$ killed at time $\boldsymbol{e},$ via the invariant
function $h^{\scriptscriptstyle{\downarrow}}.$ Under $\prc$ the
canonical process is still a L\'evy process with infinite lifetime
that drifts to $-\infty.$ Furthermore, by the third hypothesis in
(H") we have that $m^{\downarrow}=\erc(\xi_1)\in]-\infty,0[.$ We are
interested in the pssMp $(X,\pc),$ which is the Markov process
associated to the L\'evy process with law $\prc$ via Lamperti's
transformation. Since the L\'evy process $(\xi,\prc)$ drifts to
$-\infty$ we have that $(X,\pc_{x})$ hits 0 continuously at some
finite time a.s. for every $x>0.$ As a consequence of the following
result we will refer to $(X,\pc)$ as the process $(X,\p)$
conditioned to hit $0$ continuously.

\begin{teo}\label{prop:condch4} Assume that the hypotheses $(H")$ are
satisfied.
\begin{enumerate}
\item[$(i)$] For every $x>0,$  $\pcx$ is the unique measure such that
for every stopping time $T$ of $(\mathcal{G}_t)$ we have
$$\pcx(F_T, T<T_0)=x^{-\gamma}\p_x(F_T X^{\gamma}_T,T<T_0),$$ for
every $F_T\in\mathcal{G}_T.$ \item[$(ii)$] For every $x>0,$
$$\lim_{\epsilon\to 0}\p_x(F_t\cap\{t<T_0\} | X_{T_0-}\leq
\epsilon)=\pcx(F),\qquad F_t\in\mathcal{G}_{t}, \ t\geq 0.$$
\item[$(iii)$] For every $x>0,$
$$\lim_{\epsilon\to 0}\p_x(F_t\cap\{t<T_0\}  | \inf_{0\leq t<T_{0}}
X_t < \epsilon)=\pcx(F),\qquad F_t\in\mathcal{G}_{t}, \ t\geq 0.$$
\end{enumerate}
\end{teo}

\begin{proof} Part $(i)$  is an
immediate consequence of the fact that $\prc$ is an $h$--transform.
To prove (ii) we will need the following Lemma in which we determine
the tail distribution of a L\'evy process at given exponential time.

\begin{lemma}\label{lem:tailch4} Let $\sigma$ be a L\'evy process of law
$P,$ and with infinite lifetime. Assume that $\sigma$ is
non--arithmetic and that there exists a $\vartheta>0$ for which
$1<E(e^{\vartheta\sigma_1})<\infty,$ and
$E(\sigma^+_1e^{\vartheta\sigma_1})<\infty.$ Let $T_{\lambda}$ be an
exponential random variable of parameter $\lambda=\log
E(e^{\vartheta\sigma_1})$ and independent of $\sigma.$ We have that
$$\lim_{x\to\infty}e^{\vartheta x}P(\sigma_{T_{\lambda}}\geq
x)=\frac{\lambda}{\mu^{\natural}\vartheta},$$ with
$\mu^{\natural}=\er(\sigma_1 e^{\vartheta\sigma_1}).$ \end{lemma}
Lemma~\ref{lem:tailch4} is a consequence of the renewal theorem for
real--valued random variables and Cramer's method, see e.g.
Feller~\cite{MR42:5292}~\S XI.6.

\begin{proof} Observe that the function $Z(x)=P(\sigma_{T_{\lambda}}\geq x),$
satisfies a renewal equation. More precisely, for $z(x)=\int^1_0dt
\lambda e^{-\lambda t}P(\sigma_t\geq x)$ and
$L(dy)=e^{-\lambda}P(\sigma_1\in dy)$ we have that
$$Z(x)=z(x)+\int^{\infty}_{-\infty}L(dy)Z(x-y).$$ This is an elementary
consequence of the fact that the process $(\sigma'_s=
\sigma_{1+s}-\sigma_1, s\geq 0)$ is a L\'evy process independent of
$(\sigma_r,r\leq 1)$ with the same law as $\sigma.$
Next, the measure $L$ is a defective law, $L(\re)<1$,  such that
$$\int^{\infty}_{-\infty} e^{\vartheta y} L(dy) =
e^{-\lambda}E(e^{\vartheta\sigma_1})=1;\quad\text{and}\quad
\int^{\infty}_{-\infty}
y e^{\vartheta y} L(dy) < \infty,$$ by hypotheses. Thus the function
$Z^{\natural}(x)\equiv e^{\theta x}Z(x), x\in\re$ satisfies a
renewal equation with $L(\mathrm{d}y)$ replaced by
$L^{\natural}(\mathrm{d}y)=e^{\theta y}L(\mathrm{d}y),$ $y\in\re,$
and $z$ replaced by $z^{\natural}(x)=e^{\theta x}z(x),$ $x\in\re.$
By the uniqueness of the solution of the renewal equation we have
that
$$Z^{\natural}(y)=\int_{\re}z^{\natural}(y-x)U^{\natural}(\mathrm{d}x),\qquad
y\in\re,$$ where $U^{\natural}(\mathrm{d}x)$ is the renewal
measure associated to the law $L^{\natural}.$ Furthermore, the
function $z^{\natural}$ is directly Riemann integrable because it
is the product of an exponential function and a decreasing one and
$z^{\natural}$ is integrable. To see that $z^{\natural}$ is
integrable, use the Fubini's theorem to establish
\begin{equation*}
\begin{split}
\int^{\infty}_{-\infty}z^{\natural}(x)dx&=\int^1_0dt\lambda
e^{-\lambda t}E\left(\int^{\infty}_{-\infty}dx e^{\vartheta x}
1_{\{\sigma_t\geq x\}}\right)\\
&=\frac{1}{\vartheta}\int^1_0 dt \lambda e^{-\lambda
t}E(e^{\vartheta \sigma_t})\\ &=\frac{\lambda}{\vartheta}<\infty.
\end{split}
\end{equation*}
Finally, given that $L^{\natural}$ is a non-defective law with
finite mean the Key renewal theorem implies that
\begin{equation*}
\lim_{y\to\infty}Z^{\natural}(y)=\lim_{y\to\infty}
\int_{\re}z^{\natural}(y-x)U^{\natural}(\mathrm{d}x)=
\frac{1}{\mu^{\natural}}\int^{\infty}_{-\infty}z^{\natural}(x)dx=
\frac{\lambda}{\vartheta\mu^{\natural}}.
\end{equation*}\end{proof}
\noindent Now we may end the proof of part $(ii)$. Observe that
under $\p_x$ the random variable $X_{T_0-}$ has the same law as
$xe^{\xi_{\boldsymbol{e}}}$ under $\pr$. Then, applying
Lemma~\ref{lem:tailch4} to $(-\xi,\pr)$ we obtain by hypotheses~(H")
that
$$\lim_{y\to\infty}e^{-\gamma y}\pr(\xi_{\boldsymbol{e}}\leq
-y)=\frac{\ke}{\gamma\mu^{\scriptscriptstyle{\downarrow}}}:=\mathrm{d}_{\ke},$$
with
$\mu^{\scriptscriptstyle{\downarrow}}=\er(\xi_1e^{\gamma\xi_1})\in]-\infty,0[,$
which is finite by hypothesis.  Thus, we have the following estimate
of the left tail distribution of $X_{T_0-}$
\begin{equation}\label{tailX:ch4}\lim_{\epsilon\to
0}\epsilon^{\gamma}\p_x(X_{T_0-}\leq\epsilon)=x^{\gamma}\mathrm{d}_{\ke}.
\end{equation}
We conclude by a standard application of the Markov property,
estimate~(\ref{tailX:ch4}) and a dominated convergence argument.

Now we prove part $(iii)$. First of all, we claim that under the
assumptions of Theorem~\ref{prop:condch4},
\begin{equation}\label{mintailestimate}x^{-\gamma}\lim_{\epsilon\to 0+}
\epsilon^{\gamma}\p_x(\inf_{0\leq t <T_0}X_t
<\epsilon):=d^{"}_{\ke}\in]0,\infty[,\qquad x>0.\end{equation}Owing
to this estimate the rest of the proof of
Theorem~\ref{prop:condch4}~(iii) is quite similar to the one of (ii)
in Theorem~\ref{convmin0} in the case where Cramer's condition is
satisfied, so we omit the details. Indeed, it is clear that the r.v.
$\inf_{0\leq t <T_0}X_t,$ has the same law as
$$\exp\{\inf_{0\leq s\leq \boldsymbol{e}}\{\xi_{s}\}\},$$ under
$\pr.$ Its well known that $(\sup_{0\leq s\leq
\boldsymbol{e}}\{-\xi_{s}\},\pr)$ has the same law as a
subordinator, say $\widetilde{\sigma},$ with Laplace exponent
$\widehat{\kappa}(\ke,\lambda)-\widehat{\kappa}(\ke,0),$ evaluated
at an independent exponential time of parameter
$\widehat{\kappa}(\ke,0),$ where $\widehat{\kappa}(\cdot,\cdot)$ is
the bivariate Laplace exponent of the dual ladder height process
associated to $(\xi,\pr),$ see e.g. \cite{be}~Section~VI.2. So in
order to deduce the assertion~(\ref{mintailestimate}) using
Lemma~\ref{lem:tailch4} we have to verify that
\begin{equation*}
\text{(a)}~1<\er(e^{\widehat{\gamma}\widetilde{\sigma}_1})<\infty,\quad
\text{(b)}~\er(\widetilde{\sigma}_1e^{\widehat{\gamma}
\widetilde{\sigma}_1})<\infty\quad\text{and}\quad
\text{(c)}~\widehat{\kappa}(\ke,0)=\log
\er(e^{\widehat{\gamma}\widetilde{\sigma}_1}),\quad\text{for}\
\widehat{\gamma}=-\gamma.
\end{equation*}
Recall that a function $f:\re\to\re,$ of the type
$f(x)=|x|^{a}e^{\beta x},$ for $a\in\re,$ $\beta<0,$ is integrable
w.r.t. the law of $(\xi_t,\pr)$ for some $t>0$ if and only if
$f(x)1_{\{x<-1\}}$ is integrable w.r.t. the L\'evy measure of
$(\xi,\pr),$ see e.g. \cite{sa} Proposition~25.4. Furthermore,
Vigon~\cite{thesevigon} Section 6.2, established that
$f(x)1_{\{x<-1\}},$ is integrable w.r.t. the L\'evy measure of
$(\xi,\pr)$ if and only if $f(-x)1_{\{-x>1\}},$ is integrable w.r.t.
the L\'evy measure of the dual ladder height subordinator associated
to $(\xi,\pr).$ So, that (a) and (b) are consequences of the
hypotheses~(H") and the fact that the subordinator
$\widetilde{\sigma}$ has the same L\'evy measure and drift term as
the dual ladder height subordinator associated to $(\xi,\pr).$
Finally, the assertion in (c) is an easy consequence of the
inversion theorem in Vigon~\cite{thesevigon}~Section 4.3.
\end{proof}
\end{section}

\end{document}